\documentclass{article}
\usepackage{amsmath}
\usepackage{amsfonts}
\usepackage{amssymb}
\usepackage{amsthm}

\newtheorem{proposition}{Proposition}
\newtheorem{theorem}{Theorem}

\newtheorem{remark}{Remark}
\newtheorem{lemma}{Lemma}
\newtheorem{corollary}{Corollary}
\newtheorem{hypothesis}{Hypothesis}
\DeclareMathOperator{\Rep}{Rep}

\title{Special characters on star graphs and representations of
  $*$-algebras} 
\author{Vasyl Ostrovskyi\thanks{This work was supported by the DFG,
  grant no.\ 436UKR 113/71/0-1 and State Foundation for Fundamental
  Research of Ukraine, grant no. 01.07/071}} 

\date{{\em Institute of Mathematics, National Acad. Sci. of Ukraine\\
E-mail}: \texttt{vo@imath.kiev.ua}}

\begin{document}
\maketitle

\section*{Introduction}

In a series of recent papers (see
\cite{kru_rab_sam02,kru_pop_sam03,sam_zav05mfat,vla_mel_sam05} and
references therein), 
families of self-adjoint operators $A_1$, \dots, $A_n$ were studied
such that $A_1+\dots+A_n=\lambda I$ and the spectrum of each operator
$A_l$, $l=1$, \dots, $n$,  lies in fixed finite set $M_l$. Such
families can be naturally regarded as $*$-representation of certain
algebra related to a star-shaped graph $\Gamma$ and a positive function
(character)  on the set of its vertexes.

While the cases where $\Gamma$ is a Dynkin graph or extended Dynkin
graph were studied in details (see
\cite{kru_pop_sam03,mel03,mel04,ost04umzh}  and others), in the case of
general graph only partial results are known
\cite{kru_rab_sam02,mel_rab_sam04,red_roi04}. 

In the case of extended Dynkin graphs, well-known are special
characters, which possess extra invariance properties.  In this paper,
we introduce such special characters for general star-shaped graphs
and study their properties (Sec.~3). The formulas  based on the
positive solutions of certain equation, which is related to the graph
(Sec.~2). Following \cite{red_roi04} we
decompose the special character into odd and even parts and study
their evolution under reflections (Sec.~4). Notice that similar
formulas for the evolution of the even and odd parts have been
obtained independently in \cite{red05}.  The 
obtained formulas are 
used in Sec.~5 to prove that for any graph containing an extended
Dynkin graph as a proper subgraph there exists a character such that
the corresponding algebra has infinite-dimensional irreducible
$*$-representation (the converse statement that for the case of
extended Dynkin graph all irreducible representations are
finite-dimensional, was proved in the previous paper \cite{ost05})
Notice that this fact was formulated as a hypothesis in
\cite{red_roi04}.  

\section{Preliminaries}
Let $\Gamma$ be a star-shaped graph, i.e., simply-laced non-oriented
graph consisting of  $n$  branches, such that $l$-th branch has
$k_l+1$ vertexes, $l=1$, \dots, $n$, and all branches are connected at
a single root vertex (which is extreme point for all branches). 
Any function on the set of vertexes, is defined by a \emph{character} on
the graph, 
\begin{equation}
\chi=(\alpha_1^{(1)},\dots,\alpha_{k_1}^{(1)}; \dots;
\alpha_1^{(n)},\dots,\alpha_{k_n}^{(n)}; \lambda), 
\end{equation} 
where $\alpha_1^{(l)}$, \dots, $\alpha_{k_l}^{(l)}$ are the values on the
veriest of $l$-th branch starting from non-root extreme point, and
$\lambda$ is the value on the root vertex.

Given a graph $\Gamma$ and a character $\chi$ on $\Gamma$ consider
the $*$-algebra $\mathcal A_{\Gamma,\chi}$ over $\mathbb C$ generated
by self-adjoint elements $a_l$, $l=1$, \dots, $n$, which satisfy
relations 
\[
p_l(a_l)=0, \quad \sum_{l=1}^n a_l =\lambda e,
\]
where $p_l(x) = x(x-a_1^{(l)})\dots(x-a_{k_l}^{(l)})$, $l=1$, \dots,
$n$. Then $*$-representations of this algebra are $n$-tuples $A_1$,
\dots, $A_n$ of self-adjoint operators in a Hilbert space, such that
$A_1+\dots+A_n=\lambda I$ and the spectrum of each $A_l$ is contained
in $\{a_1^{(l)},\dots,a_{k_l}^{(l)}\}$.

The essential tool to study $*$-representations of $\mathcal
A_{\Gamma,\chi}$ are reflection (Coxeter)  functors
introduced in \cite{kru02} (for non-involutive case, see
\cite{ber_gel_pon73}). Namely 
there exist two functors, $S\colon \Rep \mathcal
A_{\Gamma,\chi} \to \Rep \mathcal A_{\Gamma,\chi'}$
and $T\colon \Rep \mathcal
A_{\Gamma,\chi} \to \Rep \mathcal
A_{\Gamma,\chi''}$, where 
\begin{align*}
\chi'&=(\alpha_{k_l}^{(1)}-\alpha_{k_l-1}^{(1)},\dots,
\alpha_{k_l}^{(1)}-\alpha_{0}^{(1)}; \dots; \alpha_{k_l}^{(n)}-
\alpha_{k_l-1}^{(n)},\dots, \alpha_{k_l}^{(n)}-\alpha_{0}^{(n)};\lambda'),
\\
\lambda' &= \alpha_{k_1}^{(1)} + \dots +\alpha_{k_n}^{(n)} -\lambda,
\\
\chi''&=  (\lambda-\alpha_{k_1}^{(1)},\dots, \lambda -\alpha_1^{(1)};
\dots ;\lambda-\alpha_{k_n}^{(n)},\dots, \lambda -\alpha_1^{(n)};\lambda)
\end{align*} 
(we put $\alpha_0^{(l)}=0$, $l=1$, \dots, $n$.)
The action of these functors on $*$-representations gives rise to the
action on the set of characters, $S\colon \chi \mapsto \chi'$,
$T\colon \chi\mapsto \chi''$. This action is
extensively used to study the structure and properties of
representations of algebras introduced above.

\section{Equation related to a graph}
Let $n$ and $k_l$, $l=1$, \dots, $n$ be given natural numbers, and
let $\Gamma$ be the corresponding star-shaped graph. 
In what follows, we will use solutions of the following equation
\begin{equation}\label{graph:eq_with_t}
n-1-t = \sum_{l=1}^n \frac{1}{1+t + \dots + t^{k_l}}, \quad t\ge 0.
\end{equation}

\begin{theorem}\label{th:numbers} \textup{(See \cite{ost05}).} 
The equation \eqref{graph:eq_with_t} have no solutions on $[0,\infty)$
if and only if 
the corresponding graph $\Gamma$ is one of the Dynkin graphs $A_d$,
$d\ge1$, $D_d$, $d \ge 4$, $E_6$, $E_7$, or $E_8$.
The equation \eqref{graph:eq_with_t} has a unique solution on
$[0,\infty)$ if and only 
if the corresponding graph $\Gamma$ is one of the extended Dynkin
graphs $\tilde D_4$, $\tilde E_6$, $\tilde E_7$, or $\tilde E_8$,
this solution is $t=1$. 
In all other cases (i.e., where $\Gamma$ is neither a Dynkin graph
nor an extended Dynkin graph), the equation \eqref{graph:eq_with_t}
has on $[0,\infty)$ two solutions $0<t_1<1<t_2<n-1$.
\end{theorem}

We show that $t_1t_2=1$.

\begin{proposition}\label{prop:inverse}
Let $t$ be a solutions of
\eqref{graph:eq_with_t}. Then $t^{-1}$ is also a solution of
\eqref{graph:eq_with_t}. 
\end{proposition}

\begin{proof}
Indeed,
\begin{align*}
&t^{-1} + \sum_{l=1}^n \frac{1}{1+t^{-1}+\dots +t^{-k_l}}  = t^{-1}\Bigl(1 +
\sum_{l=1}^n \frac{t^{k_l+1}}{1+t+\dots+t^{k_l}}\Bigr)
\\
&=t^{-1}\Bigl(1+\sum_{l=1}^n \frac{t^{k_l+1}}{1+t+\dots+t^{k_l}} 
+n-1-t -\sum_{l=1}^n \frac{1}{1+t+\dots+t^{k_l}} \Bigr)
\\
&= t^{-1}\Bigl( n-t + \sum_{l=1}^n
\frac{t^{k_l+1}-1}{1+t+\dots+t^{k_l}}\Bigr) = n-1.\qedhere
\end{align*}
\end{proof} 

\begin{hypothesis}
Let $r$ be the largest eigenvalue of a star-shaped graph
$\Gamma$. Then $t+t^{-1}+2 =r^2$, where $t$ is a positive solution of
\eqref{graph:eq_with_t}. 
\end{hypothesis}

\section{Special characters}
For graphs which are not Dynkin graphs, we introduce special character
as follows:
\begin{gather}
\label{special:general}
\chi_\Gamma = (\alpha_1^{(1)}, \dots, \alpha_{k_1}^{(1)}; \dots;
\alpha_1^{(n)},
\dots,\alpha_{k_1}^{(n)};1),\notag
\\
\alpha_j^{(l)}=\frac{1+t+\dots+t^{j-1}}{1+t+\dots+t^{k_l}},
\end{gather}
where $t$ is a solution of \eqref{graph:eq_with_t}. For extended
Dynkin graph such character is unique (since  \eqref{graph:eq_with_t}
has a unique solution  $t=1$) and has the form 
\begin{align*}
\chi_{\tilde D_4}& =\frac12 (1;1;1;1;2);
\\
\chi_{\tilde E_6}&=\frac13(1,2;1,2;1,2;3);
\\
\chi_{\tilde E_7}& = \frac14(1,2,3;1,2,3;2;4);
\\
\chi_{\tilde E_8}& = \frac16(1,2,3,4,5;2,4;3;6).
\end{align*}
 For any graph containing extended Dynkin graph there are two special
 characters corresponding to two positive solutions of
 \eqref{graph:eq_with_t}. 

\begin{proposition}
$(TS)\chi_\Gamma= t^{-1} \chi_\Gamma$. 
\end{proposition}

\begin{proof}
Indeed, the equality
\begin{equation}
\sum_{l=1}^n \frac{t^{k_l}}{1+\dots+t^{k_l}} = n-1-\frac1t \notag
\end{equation}
implies that
\begin{equation}\label{eq:sum_of_highest}
\sum_{l=1}^n \frac{1+\dots+t^{k_l-1}}{1+\dots+t^{k_l}} = 1+\frac1t,
\end{equation}
and by a direct calculations we get
$(TS)\chi_\Gamma=\frac1t\chi_\Gamma$.
\end{proof}

Given a graph $\Gamma$ consider the following quadratic form on the
set of characters~$\chi$:
\begin{equation}\label{eq:form}
\gamma_\Gamma(\chi)= \sum_{l=1}^n\sum_{j=1}^{k_l} (\alpha_j^{(l)})^2
+\lambda^2 - \sum_{\substack{l=1\\k_l>1}}^n 
\sum_{j=1}^{k_l-1} x_j^{(l)}x_{j+1}^{(l)} -\lambda \sum_{l=1}^{n}
x_{k_l}^{(l)}  
\end{equation}

\begin{proposition}
For special character, $\gamma_\Gamma (\chi_\Gamma) =0$.
\end{proposition}

\begin{proof}
We have, using \eqref{eq:sum_of_highest}
\begin{align*}
\gamma_\Gamma(\chi_\Gamma,1)& = 
\sum_{l=1}^n \sum_{j=1}^{k_l}\frac{(1+\dots +t^{j-1})^2} {(1+\dots
  +t^{k_l})^2} +1 
\\
&{}\qquad-\sum_{\substack{l=1\\k_l>1}}^n \sum_{j=1}^{k_l-1}\frac{(1+\dots
  +t^{j-1}) (1+\dots  +t^{j})}{(1+\dots +t^{k_l})^2}  - \sum_{l=1}^n
\frac{1+\dots +t^{k_l-1}}{1+\dots +t^{k_l}}
\\
&= \sum_{l=1}^n \frac{1}{(1-t^{k_l+1})^2}\Bigl( \sum_{j=1}^{k_l}
(1-t^{j})^2 - \sum_{j=1}^{k_l-1} (1-t^{j}) (1-t^{j+1}) \Bigr)-\frac1t.  
\end{align*}
Here we assume $t\ne1$, since for $t=1$ ($\Gamma$ being an extended
Dynkin graph) the statement is well known.

For $k_l>1$ we have
\begin{align*}
 \frac{1}{(1-t^{k_l+1})^2}\Bigl( \sum_{j=1}^{k_l}
(1-t^{j})^2 - \sum_{j=1}^{k_l-1} (1-t^{j}) (1-t^{j+1}) \Bigr)=
\frac{1-t^{k_l}}{(1-t^{k_l+1})(1+t)}. 
\end{align*}
For $k_l=1$ the same hods assuming the second summand is zero. Then by
\eqref{eq:sum_of_highest} we have 
\[
\gamma_\Gamma(\chi_\Gamma)= \frac1{t+1} \sum_{l=1}^n
\frac{1-t^{k_l}}{1-t^{k_l+1}} -\frac1t =0. \qedhere
\]
\end{proof}

\section{Evolution of odd and even components}
 
In this section, similarly to \cite{red_roi04} we decompose the
special character into odd and even parts (in \cite{red_roi04} these
parts are called odd and even standard characters) and study their evolution
under the action of the Coxeter functors. We obtain for arbitrary
star-shaped graphs the formulas
obtained in \cite{red_roi04} for special classes of graphs. 

As shown in \cite{kru_pop_sam03}, there is a correspondence between
non-degenerate $*$-represen\-ta\-tions of algebras $\mathcal
A_{\Gamma,\chi}$ and non-degenerate locally scalar representations of~
$\Gamma$. Namely, $*$-representations of the algebra $\mathcal
A_{\Gamma,\chi}$ correspond to locally scalar representations of the
graph $\Gamma$ with character
\begin{equation}\label{eq:standard:ls0}
u_\Gamma=(x_1^{(1)},\dots,x_{k_1}^{(1)};
\dots;x_1^{(n)},\dots,x_{k_n}^{(n)};1)  
\end{equation}
where
\begin{align}\label{chi2u}
x_{k_l}^{(l)}& = \alpha_{k_l}^{(l)}, &
x_{k_l-1}^{(l)}& = \alpha_{k_l}^{(l)}-\alpha_{1}^{(l)}, \notag
\\
x_{k_l-2}^{(l)}& = \alpha_{k_l-1}^{(l)}-\alpha_{1}^{(l)}, 
&
x_{k_l-3}^{(l)}& = \alpha_{k_l-1}^{(l)}-\alpha_{2}^{(l)}, \notag
\\
x_{k_l-4}^{(l)}& = \alpha_{k_l-2}^{(l)}-\alpha_{2}^{(l)}, 
&
x_{k_l-5}^{(l)}& = \alpha_{k_l-2}^{(l)}-\alpha_{3}^{(l)}, 
\end{align}
and so on; $l=1$, \dots, $n$.

Applying these formulas for the special character $\chi_\Gamma$ we
have the following special character of locally scalar representation: 
\begin{align}
x_{k_l}^{(l)}& = 
\frac{1+t+\dots+t^{k_l-1}} {1+t+\dots+ t^{k_l}},
&
x_{k_l-1}^{(l)}&=
\frac{t\,(1+t+\dots+t^{k_l-2})} {1+t+\dots+ t^{k_l}},
\notag
\\
x_{k_l-2}^{(l)}&=
\frac{t\,(1+t+\dots+t^{k_l-3})} {1+t+\dots+ t^{k_l}}.
&
x_{k_l-3}^{(l)}&=
\frac{t^2\,(1+t+\dots+t^{k_l-4})} {1+t+\dots+ t^{k_l}},
\label{eq:standard:ls}
\end{align}
and so on. 

Following \cite{kru_roi05} we decompose the set of vertexes
$\Gamma_v$ into odd and even parts,
$\Gamma_v=\overset{\bullet}\Gamma_v\cup
\overset{\circ}\Gamma_v$. However, unlike in \cite{kru_roi05},  
we always  assume that the root vertex is odd. According to the
decomposition of $\Gamma_v$ we decompose $u_\Gamma$ into odd end even
parts, $u_\Gamma=\overset{\bullet}u_\Gamma +\overset{\circ}u_\Gamma$,
where 
\begin{align}
\overset{\bullet}u_\Gamma & =
(\dots,0,x_{k_1-3}^{(1)},0,x_{k_1-1}^{(1)},0; \dots;
\dots,0,x_{k_n-3}^{(n)},0,x_{k_n-1}^{(n)},0;1),\notag 
\\
\overset{\circ}u_\Gamma & =
(\dots,0,x_{k_1-2}^{(1)},0,x_{k_1}^{(1)}; \dots;
\dots,0,x_{k_n-2}^{(n)},0,x_{k_n}^{(n)};0). \notag
\end{align}

Denote by $\sigma_g$ the reflection at a point $g\in \Gamma_v$. 
Let $\overset{\bullet}c$ and $\overset{\circ}c$ be compositions of
reflections at all odd and even points respectively (the order is
irrelevant since all odd reflections commute and all even reflections
commute). For detailed definitions and properties of these mappings,
see \cite{kru_roi05}.

\begin{proposition}\label{prop:sigma}
Let the character of locally scalar representation be defined by
\eqref{eq:standard:ls0}, \eqref{eq:standard:ls}. For odd vertexes
$\sigma_g(x_g) = t^{-1} 
x_g$. For even vertexes $\sigma_g(x_g) = t x_g$.
\end{proposition}

\begin{proof}
Straightforward calculations.
\end{proof}

\begin{proposition}\label{prop:cox_standard}
The following formulas hold for all $j=1$, $2$,~\dots.
\begin{align*}
(\overset{\bullet}c\,\overset{\circ}c)^j\overset{\bullet}u_\Gamma
&= 
\frac{1}{t^j(1-t)} \bigl(
(1-t^{2j+1})\overset{\bullet}u_\Gamma +
t (1-t^{2j})\overset{\circ}u_\Gamma \bigr), 
\\
\overset{\circ}c
(\overset{\bullet}c\,\overset{\circ}c)^j\overset{\bullet}u_\Gamma
&= 
\frac{1}{t^{j}(1-t)} \bigl(
(1-t^{2j+1})\overset{\bullet}u_\Gamma +
 (1-t^{2j+2})\overset{\circ}u_\Gamma \bigr), 
\\
(\overset{\circ}c\,\overset{\bullet}c)^j\overset{\circ}u_\Gamma
&= 
\frac{1}{t^{j}(1-t)} \bigl(
(1-t^{2j})\overset{\bullet}u_\Gamma +
 (1-t^{2j+1})\overset{\circ}u_\Gamma \bigr), 
\\
\overset{\bullet}c
(\overset{\circ}c\,\overset{\bullet}c)^j\overset{\circ}u_\Gamma
&= 
\frac{1}{t^{j+1}(1-t)} \bigl(
(1-t^{2j+2})\overset{\bullet}u_\Gamma +
 t(1-t^{2j+1})\overset{\circ}u_\Gamma \bigr). 
\end{align*}
\end{proposition}

\begin{proof}
The  following fact follows from the Proposition~\ref{prop:sigma}
above. 

\begin{lemma}\label{prop:cd_cc}
Let $u_\Gamma$ be defined by \eqref{eq:standard:ls0},
\eqref{eq:standard:ls}. Then 
$\overset{\bullet}{c} \overset{\bullet}u_\Gamma
=-\overset\bullet u_\Gamma$, 
$\overset{\circ}{c} \overset{\circ}u_\Gamma
=-\overset{\circ} u_\Gamma$, 
$\overset{\circ}{c} \overset{\bullet}u_\Gamma
=\overset{\bullet}u_\Gamma +(1+t)\overset{\circ}u_\Gamma$, 
$\overset{\bullet}{c} \overset{\circ} u_\Gamma
=\overset\circ u_\Gamma +(1+t^{-1})\overset{\bullet} u_\Gamma$, 
\end{lemma}

Now the statement follows by induction from Lemma~\ref{prop:cd_cc}.
\end{proof}

\begin{remark}
Since $\overset{\bullet}c \overset{\bullet} u_\Gamma = -
\overset{\bullet} u_\Gamma$ and $\overset{\circ}c
\overset{\circ} u_\Gamma =- \overset{\circ} u_\Gamma$,
Proposition~\ref{prop:cox_standard} obviously gives the formulas for
$(\overset{\circ}c\,\overset{\bullet}c)^j
\overset{\bullet} u_\Gamma$,
$\overset{\bullet}c(\overset{\circ}c\,\overset{\bullet}c)^j
\overset{\bullet} u_\Gamma$,
$(\overset{\bullet}c\,\overset{\circ}c)^j
\overset{\circ} u_\Gamma$,
$\overset{\circ}c(\overset{\bullet}c\,\overset{\circ}c)^j
\overset{\circ} u_\Gamma$ as well. 
\end{remark}

\begin{corollary}
For a character of the form $v=\alpha\overset{\bullet}u_\Gamma +
\beta\overset{\circ}u_\Gamma$ ($\alpha\ne0$)
introduce a ``normalized'' character $\tilde
v=\overset{\bullet}u_\Gamma + 
\alpha^{-1}\beta\overset{\circ}u_\Gamma$. Then we have formulas 
similar to \cite{red_roi04}:
\begin{align*}
\widetilde{(\overset{\bullet}c\,\overset{\circ}c)^j
  \overset{\bullet}u_\Gamma}  
&
=  \overset{\bullet}u_\Gamma +
\frac{t}{1+t} \rho_{\lambda-2}(2j)\,\overset{\circ}u_\Gamma,
\\
\widetilde{\overset{\circ}c
(\overset{\bullet}c\,\overset{\circ}c)^j \overset{\bullet}u_\Gamma} 
&
=\overset{\bullet}u_\Gamma +
(1+t)(\rho_{\lambda-2}(2j+1))^{-1}\,\overset{\circ}u_\Gamma
\\
&=\overset{\bullet}u_\Gamma +
\frac{t}{1+t}(\lambda+2-\rho_{\lambda-2}(2j))\,
  \overset{\circ}u_\Gamma ,
\\
\widetilde{(\overset{\circ}c\,\overset{\bullet}c)^j
  \overset{\circ}u_\Gamma} 
&
=\overset{\bullet}u_\Gamma +
(1+t)(\rho_{\lambda-2}(2j))^{-1}\overset{\circ}u_\Gamma
\\
&=\overset{\bullet}u_\Gamma +
\frac{t}{1+t}(\lambda+2-\rho_{\lambda-2}(2j-1))\,
  \overset{\circ}u_\Gamma ,
\\
\widetilde{\overset{\bullet}c
(\overset{\circ}c\,\overset{\bullet}c)^j \overset{\circ}u_\Gamma} 
&
=  \overset{\bullet}u_\Gamma +
\frac{t}{1+t} \rho_{\lambda-2}(2j+1)\,\overset{\circ}u_\Gamma.
\end{align*}
Here $\lambda=t+t^{-1}$ and following \cite{red_roi04} we denote 
\[
\rho_{\lambda-2}(n)=1+\frac{t-t^n}{1-t^{n+1}}. 
\]
We use the formula
\[
(\rho_{\lambda-2}(n+1))^{-1} = \frac{t}{(t+1)^2} (\lambda +2
-\rho_{\lambda-2}(n)), 
\]
which is verified directly.
Notice that in \cite{red_roi04}, the condition $t>1$ was assumed,
while we admit as $t$ any of the two positive solutions of
\eqref{graph:eq_with_t}.  
\end{corollary}

\begin{proposition}\label{prop:limit_u} Let $u_\Gamma'$ be the
  special character, 
  constructed by the same formulas \eqref{eq:standard:ls0},
  \eqref{eq:standard:ls} as $u_\Gamma$ 
  but with $t^{-1}$ instead of $t$. Then for $t>1$
\begin{align*}
\lim_{j\to \infty} \widetilde{(\overset{\bullet}c\,\overset{\circ}c)^j
  \overset{\bullet}u_\Gamma} 
& =\lim_{j\to\infty}\widetilde{\overset{\bullet}c 
(\overset{\circ}c\,\overset{\bullet}c)^j \overset{\circ}u_\Gamma}
=u_\Gamma, 
 \\
\lim_{j\to\infty} \widetilde{\overset{\circ}c
(\overset{\bullet}c\,\overset{\circ}c)^j
  \overset{\bullet}u_\Gamma}
& =\lim_{j\to\infty} \widetilde{(\overset{\circ}c\,\overset{\bullet}c)^j
  \overset{\circ}u_\Gamma} =u_\Gamma',
\end{align*}
and for $t<1$
\begin{align*}
\lim_{j\to \infty} \widetilde{(\overset{\bullet}c\,\overset{\circ}c)^j
  \overset{\bullet}u_\Gamma} 
& =\lim_{j\to\infty}\widetilde{\overset{\bullet}c 
(\overset{\circ}c\,\overset{\bullet}c)^j \overset{\circ}u_\Gamma}
=u_\Gamma', 
 \\
\lim_{j\to\infty} \widetilde{\overset{\circ}c
(\overset{\bullet}c\,\overset{\circ}c)^j
  \overset{\bullet}u_\Gamma}
& =\lim_{j\to\infty} \widetilde{(\overset{\circ}c\,\overset{\bullet}c)^j
  \overset{\circ}u_\Gamma} =u_\Gamma.
\end{align*}
\end{proposition}

\begin{proof}
Let $t>1$. It was shown in \cite{red_roi04} that $\lim_{j\to\infty}\rho_{\lambda-2}(j)
  =1 +t^{-1}$ which implies the first formula.

The same formula yields 
\[
\lim_{j\to\infty} \widetilde{\overset{\circ}c
(\overset{\bullet}c\,\overset{\circ}c)^j
  \overset{\bullet} u_\Gamma} =\overset{\bullet}u_\Gamma +
t\, \overset{\circ} u_\Gamma
\]

Let $\smash{\overset{\bullet}u}_\Gamma'$ and
$\smash{\overset{\circ}u}_\Gamma'$
be odd and even components of $u_\Gamma'$. The nonzero components
of  $\overset{\bullet}u_\Gamma$ have the form
\[
x_{k_l-2j+1}^{(l)} =\frac{t^j(1+t+\dots+t^{k_l-2j})}{1+t+\dots
  +t^{k_l}}, \quad
j=0,1,\dots 
\]
which is invariant under replacement $t\mapsto t^{-1}$. Therefore,
$\smash{\overset{\bullet}u}_\Gamma' = \overset{\bullet}u_\Gamma$. 

The nonzero components of $\overset{\circ}u_\Gamma$ have the form
\[
x_{k_l-2j}^{(l)} = \frac{t^j(1+t+\dots+t^{k_l-2j-1})}{1+t+\dots
  +t^{k_l}},\quad j=0,1,\dots 
\]
which implies $\smash{\overset{\circ}u}'_\Gamma  =t\,
\overset{\circ}u_\Gamma$.

For $t<1$ we have $\lim_{j\to\infty}\rho_{\lambda-2}(j)
  =1 +t$. The rest of the proof is the same.
\end{proof}

\section{Representations}

\begin{proposition}
Let $\Gamma$ be a star-shaped graph. If\/  $\Gamma$ contains an
extended Dynkin graph as a proper subgraph, then there exists a
character $\chi$ on $\Gamma$ such that $\mathcal A_{\Gamma,\chi}$ has 
infinite-dimensional irreducible $*$-representation. 

Otherwise, all irreducible $*$-representations of  $\mathcal
A_{\Gamma,\chi}$ are finite-dimensional regardless of the choice of
$\chi$. 
\end{proposition}

\begin{proof}
It was shown in previous paper \cite{ost05} that for Dynkin graphs and
extended Dynkin graphs all irreducible
$*$-representations of  $\mathcal A_{\Gamma,\chi}$ are
finite-dimensional. Therefore, we need to prove the first part of the
proposition only.

Let $\Gamma'\supset\Gamma$. Then any irreducible representation of  $\mathcal
A_{\Gamma,\chi}$ can trivially be extended to irreducible
representation of  $\mathcal A_{\Gamma',\chi'}$, where $\chi'$
coincides with $\chi$ on vertexes of $\Gamma$ and is arbitrary outside
$\Gamma$. Therefore, it remains to prove the statement for the
following graphs: $\Gamma_{(1,1,1,1,1)}$, $\Gamma_{(1,1,1,2)}$,
$\Gamma_{(2,2,3)}$, $\Gamma_{(1,3,4)}$, $\Gamma_{(1,2,6)}$, since any
star-shaped graph containing an extended Dynkin graph as a proper
subgraph, contains (or coinsides with) one of the five listed above.
Here we denote by $\Gamma_{(k_1,\dots,k_n)}$ the graph with $n$
branches of lengths $k_1$, \dots, $k_n$.

After some calculations one can see that the components of the
characters  
\begin{align*}
w_p=\overset{\bullet}u_\Gamma +
\frac{1-t^{p}}{1-t^{p-1}} \,\overset{\circ}u_\Gamma,
\end{align*}
are increasing along branches, 
$(w_p)_{j+1}^{(l)}>(w_p)_j^{(l)}$,    for $p\ge k=\max k_l$.
Let $\chi_p$ be the corresponding character on the
graph constructed by the procedure inverse to \eqref{chi2u}. 
Then the components of  $\chi_p$ are positive and
increasing along branches, $(\chi_p)_{j+1}^{(l)}>(\chi_p)_j^{(l)}$,
$j=1$, \dots, $k_l-1$, $l=1$, \dots, $n$.  This property of characters
is essential  for the constructions related to  the
corresponding algebra $\mathcal A_{\Gamma,\chi_p}$. 
Moreover, for $p=k$, we have $(\chi_k)_{k}^{(l)}=1$ ($k_l=k$) which
allows to construct for this algebra one-dimensional
representation.  

Since the action of $T$ and $S$  is induced by the action of
$\overset{\circ}c$ and $\overset{\bullet}c$, applying powers of
$(TS)$, to this representation, we get irreducible representations of
the whole family of 
algebras $\mathcal A_{\Gamma,\chi_j}$, $j=p+2r$, $r=0$, 1,\dots. 

Applying Proposition~\ref{prop:limit_u} and results of \cite{shu02} we 
conclude that the algebra $\mathcal A_{\Gamma, \chi_\Gamma}$ has at
  least one $*$-representation.   

If $\mathcal A_{\Gamma,\chi_\Gamma}$ has a finite-dimensional
irreducible representation, then the trace equality implies that 1
belongs to the rational span of  components of $\chi_\Gamma$.
To complete the proof, we show that this is not true for the five
graphs listed above. This is known for $\Gamma_{(1,1,1,1,1)}$
\cite{kru_rab_sam02} and for $\Gamma_{(1,1,1,2)}$ \cite{mel_rab_sam04}

Indeed, after simple transformations we find that $t$ is a root of the
corresponding polynomial $p_\Gamma(\cdot)$:
\begin{align*}
p_{\Gamma_{(1,1,1,1,1)}}(x)&=x^2 -3x +1,
\\
p_{\Gamma_{(1,1,1,2)}}(x)&=x^4 -x^3 -x^2 -x +1,
\\
p_{\Gamma_{(2,2,3)}}(x)&=x^6 -x^4 -x^3 -x^2  +1 ,
\\
p_{\Gamma_{(1,3,4)}}(x)&=x^8 -x^5 -x^4 -x^3 +1 ,
\\
p_{\Gamma_{(1,2,6)}}(x)&=x^{10} +x^9 -x^7 -x^6 -x^5 -x^4 -x^3 +x +1.
\end{align*}
It is a routine exercise to check that each of these polynomials is
irreducible over $\mathbb N$ and therefore over $\mathbb Q$. Therefore,
if $t$ is a root of 
$p_\Gamma(\cdot)$, the numbers $1$, $t$, \dots, $t^{m-1}$, $m=\deg
p_\Gamma(\cdot)$,  are rational independent.

Reducing the components of $\chi_\Gamma$ to common denominator, we
obtain the following problem: to check that

---for  $\Gamma_{(1,1,1,1,1)}$ $q_\Gamma(t)=1+t$ is rational
   independent of $1$; 

---for $\Gamma_{(1,1,1,2)}$ $q_\Gamma(t)=(1+t)(1+t+t^2)$ is rational
   independent of $1+t+t^2$, $1+t$, $(1+t)(1+t)$;

---for $\Gamma_{(2,2,3)}$ $q_\Gamma(t)=(1+t+t^2)(1+t+t^3)$ is rationally
   independent of $1+t+t^2+t^3$, $(1+t)(1+t+t^2+t^3)$, $1+t+t^2$,
   $(1+t)(1+t+t^2)$, $(1+t+t^2)(1+t+t^2)$;

---for $\Gamma_{(1,3,4)}$
   $q_\Gamma(t)=(1+t)(1+t+t^2+t^3)(1+t+t^2+t^3+t^4)$ is 
   rationally independent of $(1+t+t^2+t^3)(1+t+t^2+t^3+t^4)$,
   $(1+t)(1+t+t^2+t^3+t^4)$, $(1+t)(1+t)(1+t+t^2+t^3+t^4)$,
   $(1+t+t^2)(1+t)(1+t+t^2+t^3+t^4)$, $(1+t)(1+t+t^2+t^3+t^4)$,
   $(1+t)(1+t)(1+t+t^2+t^3+t^4)$, $(1+t+t^2)(1+t)(1+t+t^2+t^3+t^4)$,
   $(1+t+t^3+t^4)(1+t)(1+t+t^2+t^3+t^4)$;

---for $\Gamma_{(1,2,6)}$ $q_\Gamma(t)=
   (1+t)(1+t+t^2)(1+t+t^2+t^3+t^4+t^5+t^6)$ is 
   rationally independent of  $(1+t+t^2)(1+t+t^2+t^3+t^4+t^5+t^6)$
   $(1+t)(1+t+t^2)(1+t+t^2+t^3+t^4+t^5+t^6)$,
   $(1+t)(1+t+t^2+t^3+t^4+t^5+t^6)$,
   $(1+t)(1+t)(1+t+t^2+t^3+t^4+t^5+t^6)$, $(1+t)(1+t+t^2)$,
   $(1+t)(1+t)(1+t+t^2)$, $(1+t+t^2)(1+t)(1+t+t^2)$
   $(1+t+t^2+t^3)(1+t)(1+t+t^2)$, $(1+t+t^2+t^3+t^4)(1+t)(1+t+t^2)$,
   $(1+t+t^2+t^3+t^4 +t^5)(1+t)(1+t+t^2)$.  

But this immediately follows since none of these numbers except
contains $t^{m-1}$, $m=\deg p_\Gamma(\cdot)$.
\end{proof}

The author expresses his deep gratitude to Prof.~Yu.~S.~Samoilenko and
Prof.~S.~V.~Popovych  for
valuable discussions of the subject of this paper.


\begin{thebibliography}{10}

\bibitem{ber_gel_pon73}
I.~N. Bernstein, I.~M. Gelfand, and V.~A. Ponomarev, \emph{Coxeter functors and
  {G}abriel's theorem}, Uspekhi Mat. Nauk \textbf{28} (1973), no.~2, Russian.

\bibitem{kru_pop_sam03}
S.~Krugljak, S.~Popovych, and Yu. Samoilenko, \emph{Representations of
  $*$\nobreakdash-algebras associated with {D}ynkin graphs and
  {H}orn's problem}, Uchen. 
  Zap. Tavrich. Nat. Univ. \textbf{16} (2003), no.~2, 132--139.

\bibitem{kru_roi05}
S.~A. Krugljak and A.~V. Roiter, \emph{Locally scalar representations of graphs
  in the category of {H}ilbert spaces}, Funct. Anal. Prilozh. \textbf{39}
  (2005), no.~2, 13--30.

\bibitem{kru_rab_sam02}
S. A. Krugljak, V. I. Rabanovich, and Yu. S. Samoilenko, \emph{On sum of
  projections}, Funct. Anal. Prilozh. \textbf{36} (2002), no.~3, 20--35.

\bibitem{kru02}
S.~A. Kruglyak, \emph{Coxeter functors for a certain class of
  $*$\nobreakdash-quivers and 
  $*$\nobreakdash-algebras}, Methods Funct. Anal. Topol. \textbf{8}
(2002), no.~4, 49--57. 

\bibitem{mel04}
A.~Mellit, \emph{Algebras generated by elements with given spectrum and scalar
  sum and kleinian singularities}, arXiv:math.RA/0406119.

\bibitem{mel03}
A.S. Mellit, \emph{When a sum of three partial reflections is zero}, Ukr. Math.
  J. \textbf{55} (2003), no.~9, 1277--1283.

\bibitem{mel_rab_sam04}
A.S. Mellit, V.I. Rabanovich, and Yu.S. Samoilenko, \emph{When a sum of partial
  reflections is a scalar operator}, Funct. Anal. Prilozh. \textbf{38} (2004),
  no.~2, 91--94.

\bibitem{ost04umzh}
V.~L. Ostrovskyi, \emph{Representations of algebra associated with {D}ynkin
  graph ${\tilde E_7}$}, Ukr. Mathem. J. \textbf{56} (2004), no.~9, to
appear. 

\bibitem{ost05}
V.~L. Ostrovskyi, \emph{On $*$-representations of a certain class of algebras
  related to a graph}, Methods Funct. Anal. Topol. \textbf{11} (2005), no.~3,
to appear (See also \texttt{arXiv:math.RT/0506505}).  

\bibitem{red05}
I.~K. Redchuk, \emph{Separating functions, spectral graph theory and locally
  scalar representations in {H}ilbert spaces}, arXiv:math.RT/0509417, 2005.

\bibitem{red_roi04}
I.~K. Redchuk and A.~V. Roiter, \emph{Singular locally scalar representations
  of quivers in {H}ilbert spaces and separating functions}, Ukr. Math. Journ.
  \textbf{56} (2004), no.~6, 947--963.

\bibitem{sam_zav05mfat}
Yu. Samoilenko and M.~Zavodovsky, \emph{Spectral theorems for
  $*$-representations of the algebras ${\mathcal{A}_{\Gamma,\chi,com}}$
  associated with {D}ynkin graphs}, Methods Funct. Anal. Topol. (2005), no.~1.

\bibitem{shu02}
V.~S. Shulman, \emph{On representations of limit relations}, Meth. Funct. Anal.
  Topol. \textbf{7} (2002), no.~4, 85--87.

\bibitem{vla_mel_sam05}
M.~S. Vlasenko, A.~S. Mellit, and Yu.~S. Samoilenko, \emph{On algebras
  generated with linearly dependent generators that have given spectra}, Funct.
  Anal. Appl. \textbf{39} (2005), no.~3, to appear.

\end{thebibliography}

\end{document}